\documentclass[12pt,a4paper] {article}
\usepackage[cp1251]{inputenc}
\usepackage[T2A]{fontenc}
\usepackage[english]{babel}
\usepackage{amsmath}
\usepackage{amssymb}

\usepackage{geometry}
\newcommand{\bs}{\boldsymbol}
\newcommand{\ds}{\displaystyle}

\begin{document}

\begin{center}

{\bf \large An Explicit Integration of a Problem of Motion \\
of a Generalized Kovalevskaya Top\footnote{Doklady Mathematics, Vol. 71, No. 2, 2005, pp. 298–299. Translated from Doklady Akademii Nauk, Vol. 401, No. 3, 2005, pp. 321–323. }}

\vspace{2mm}
{\bf M. P. Kharlamov and A. Yu. Savushkin\footnote{Volgograd Academy of Public Administration, Chair of Computer Systems and Mathematical Simulation}}
\vspace{3mm}

Presented by Academician V.V. Rumyantsev June 18, 2004 \\
Received July 8, 2004

\end{center}

\vspace{3mm}

Kovalevskaya’s solution [1] of the problem of
motion of a heavy rigid body about a fixed point was
generalized to the case of a double constant field in [2,
3]. The corresponding Hamiltonian system has three
degrees of freedom. The invariant four-dimensional
submanifolds of the phase space were found in [2, 4].
The case of [2] was studied in [5]. In this paper, we consider
the case of [4].

Consider a rigid body having a fixed point and satisfying
the Kovalevskaya condition for the principal moments of inertia at the fixed point ${\bf
I}={\rm diag}\,\{2,2,1\}$. Suppose that the body is placed in a force field with the potential
$$
 U=-({\bf e}_1,{\boldsymbol \alpha})-({\bf
e}_2,{\boldsymbol\beta}),\eqno(1) $$
where ${\bf e}_1$ and ${\bf e}_2$ are the orthonormal vectors directed
along the principal axes of inertia in the equatorial
plane and ${\boldsymbol \alpha}$ and ${\boldsymbol \beta}$ are fixed vectors in the inertial space. 

\textbf{Theorem 1.} \textit{Without loss of generality we can assume the vectors ${\bs
\alpha},{\bs \beta}$ to be mutually orthogonal.}

{\bf Proof}. Note that the potential (1) is invariant with respect to the
substitution
$$
 \begin{pmatrix}
   {\bf e_1 }  \\
   {\bf e_2 }
 \end{pmatrix}  \mapsto \Theta
 \begin{pmatrix}
   {\bf e_1 }  \\
   {\bf e_2 }
   \end{pmatrix},\quad
   \begin{pmatrix}
   {\boldsymbol \alpha}  \\
   {\boldsymbol \beta}
 \end{pmatrix} \mapsto \Theta
 \begin{pmatrix}
   {\boldsymbol \alpha}  \\
   {\boldsymbol \beta }
 \end{pmatrix} , \quad
 \Theta  =
 \begin{Vmatrix}
   {\cos \theta } & {\sin \theta }  \cr
   { - \sin \theta } & {\cos \theta }
 \end{Vmatrix} ,\quad
\displaystyle {\rm tg}\,2 \theta  = \frac{2{\boldsymbol \alpha}
{\boldsymbol \cdot} {\boldsymbol \beta}} {{\boldsymbol \alpha}^2 -
{\boldsymbol \beta}^2 },
$$
which leaves the pair ${\bf
e}_1,{\bf e}_2$ orthonormal in the equatorial plane of the body but makes
the vectors ${\boldsymbol \alpha},{\boldsymbol \beta}$ orthogonal to each other.

Let ${\bs \alpha}^2=a^2$ and ${\bs \beta}^2=b^2$. We consider the general
case $a > b$ and denote $p^2=a^2+b^2$, $r^2=a^2-b^2$.

The corresponding Euler\,--\,Poisson equations are Liouville\,--\,Arnold integrable due to the first integrals [2, 3]
$$
\begin{array}{l}
\ds H = \omega _1^2 + \omega _2^2 + \frac{1} {2}\omega _3^2 - (\alpha
_1 + \beta _2 ), \\[3mm]
\ds K = (\omega _1^2 - \omega _2^2 + \alpha _1
- \beta _2 )^2 +
(2\omega _1 \omega _2 + \alpha _2 + \beta _1 )^2,\\[3mm]
\ds G = \frac{1}{4}(\omega_\alpha^2 + \omega_\beta^2 ) +
\frac{1}{2}\omega_3\omega_\gamma - b^2 \alpha_1 - a^2 \beta_2,
\end{array}
$$
where
$$
\begin{array}{l}
\ds \omega_\alpha = 2\omega_1 \alpha_1 + 2\omega_2 \alpha_2 + \omega_3
\alpha_3 , \quad
\omega_\beta = 2\omega_1 \beta_1 + 2\omega_2 \beta_2 + \omega_3 \beta_3 , \\
\ds \omega_\gamma = 2\omega_1 (\alpha_2 \beta_3 - \alpha_3 \beta_2 ) +
2\omega_2 (\alpha_3 \beta_1 - \alpha_1 \beta_3 ) + \omega_3
(\alpha_1 \beta_2 - \alpha_2 \beta_1 ).
\end{array}
$$

Let us introduce the functions on the phase space
$$
\begin{array}{c}
F = (2G - p^2 H)^2 - r^4 K, \quad M =(2G - p^2 H)/r^4, \\[3mm]
L = \sqrt{2p^2 M^2+2 H M+1}.
\end{array}
$$
They are also first integrals of motion. Denote by $N$ the set of critical points of the function $F$ on the level $F = 0$.

Now we consider the new variables ($i^2=-1$):
$$
\begin{array}{c}
\begin{array}{ll}
{x_1  = (\alpha _1  - \beta _2 ) + i(\alpha _2  + \beta _1 ),} &
{x_2  = \overline {x_1 } ,} \cr {y_1  = (\alpha _1  + \beta _2 ) +
i(\alpha _2  - \beta _1 ),} & {y_2  = \overline {y_1 } ,}\cr
{z_1=\alpha _3  + i\beta _3 ,} & {z_2  = \overline {z_1 } ,}
\end{array}  \\
w_1  = \omega _1  + i\omega _2, \quad w_2  = \overline {w_1},\quad
w_3  = \omega _3.
\end{array} \eqno(2)
$$

\vspace{2mm}
\textbf{Theorem 2}. {\it In the domain $x_1 x_2 \neq 0$ the set $N$ is
determined by the two independent equations
$$ F_1=0, \quad F_2=0, \eqno(3)$$
where
$$
\displaystyle {F_1 = \sqrt {x_1 x_2 } w_3 - \frac{1} {{\sqrt {x_1
x_2} }}(x_2 z_1 w_1 + x_1 z_2 w_2 ),} \quad \displaystyle {F_2 =
\frac{i} {2}[\frac{{x_2 }} {{x_1 }}(w_1^2 + x_1) - \frac{{x_1}}
{{x_2}}(w_2^2 + x_2)]}.
$$
In particular, in this domain, $N$ is a smooth four-dimensional manifold.
The induced vector field on $N$ is Hamiltonian everywhere
except at the points where $L = 0$. }

\vspace{2mm}

{\bf Proof}.
The proof of the first statement of the theorem is by
direct calculation. The second statement follows from the invariance
of the set of critical points of the first integral and
from the relation $\{F_2,F_1\}=r^2 L$ for the Poisson bracket.

\vspace{2mm}

Equations (3) were obtained in [4], but they do not
describe the invariant set $N$ as a whole because of the
presence of an obvious singularity. The theorem stated
above implies that $N$ is determined globally, and, on
this set, we have a completely integrable Hamiltonian
system with two degrees of freedom such that the set of
points where the symplectic structure degenerates is
a thin set.

It is convenient to take the functions $M$ and $L$ as an
involutive pair of first integrals.

\vspace{2mm}

{\bf Theorem 3}. {\it The change of variables
$$
s_1=\frac{x_1 x_2+z_1 z_2+r^2}{2\sqrt{x_1 x_2}},\quad
s_2=\frac{x_1 x_2+z_1 z_2-r^2}{2\sqrt{x_1 x_2}} \eqno(4)
$$
reduces the equations of motion of the Kovalevskaya
top in a double force field on the integral manifold
$$ J_{m,\ell}=\{M=m,L=\ell\}\subset N $$
to the system
$$
\begin{array}{c}
\displaystyle{\frac{ds_1}{dt}=\frac{1}{2}\sqrt{(a^2-s_1^2)\Phi(s_1)},\quad
\frac{ds_2}{dt}=\frac{1}{2}\sqrt{(b^2-s_2^2)\Phi(s_2)},}
\end{array} \eqno(5)
$$
where $\Phi(s)=4m s^2 -4 \ell s +(\ell^2-1)/m$. The solutions of this system can be written explicitly in terms of elliptic functions.}

\vspace{2mm}

{\bf Proof}.
Let us eliminate the variables $w_i$ $(i=1,2,3)$ from
relations (3) and the integral equation
$M=m$. The equation $L=\ell$ takes the form
$$
m(x_1 x_2+z_1 z_2)-\ell \sqrt{x_1 x_2}+\sqrt{ m^2 r^4 - m r^2
(x_1+x_2)+x_1 x_2}=0. \eqno(6)
$$
The vectors ${\boldsymbol \alpha},{\boldsymbol \beta}$ have constant length and are mutually orthogonal. In variables (2) this fact is written as
$$
z_1^2+x_1 y_2=r^2,\;z_2^2+x_2 y_1=r^2,\;x_1 x_2+y_1 y_2+2z_1
z_2=2p^2. \eqno(7)
$$
Then calculating the time derivatives of variables (4) and taking into account equation (6),
we obtain system (5).

\vspace{2mm}

Note that, by virtue of (7), the variables $s_1$ and $s_2$ satisfy the natural constraints
$s_1^2 \geqslant a^2,s_2^2 \leqslant b^2$.
Therefore, the real solutions of (5) oscillate in the intervals
where $\Phi(s_1)\leqslant 0, \Phi(s_2) \geqslant
0$. The separating set in the plane $(m,\ell)$ coincides with the discriminant
set of the polynomial $(a^2-s^2)(b^2-s^2)\Phi(s)$
(a system of straight lines) and a half-line $\{\ell=0,\;m <0\}$
(the latter follows directly from the definition of the function $L$).

Consider the following polynomial in two auxiliary
variables:
$$
\Psi (s_1 ,s_2 ) = 4ms_1 s_2 - 2\ell (s_1 + s_2 ) +(\ell ^2 -
1)/m.
$$
For variables (2), we obtain explicit dependencies on $s_1,s_2$:
$$
\begin{array}{l}
\displaystyle{x_1 = - \frac{{r^2 }} {{2(s_1 - s_2)^2}}[\Psi (s_1
,s_2 ) + \sqrt
{\Phi (s_1 )\Phi (s_2 )} ],}\\
\displaystyle{x_2 = - \frac{{r^2 }} {{2(s_1 -
s_2)^2}}[\Psi (s_1 ,s_2 ) - \sqrt {\Phi (s_1 )\Phi (s_2 )} ],} \\
\displaystyle{y_1 = 2\frac{{(2s_1 s_2 - p^2 ) - 2\sqrt {(s_1^2 -
a^2 )(s_2^2 - b^2 )} }} {{\Psi (s_1 ,s_2 ) - \sqrt {\Phi (s_1
)\Phi (s_2 )}
}},}\\
\displaystyle{y_2 = 2\frac{{(2s_1 s_2 - p^2 ) + 2\sqrt {(s_1^2 -
a^2 )(s_2^2 - b^2 )} }} {{\Psi (s_1 ,s_2 ) + \sqrt {\Phi (s_1
)\Phi
(s_2 )} }},} \\
\displaystyle{z_1 = \frac{r} {{s_1 - s_2 }}(\sqrt {s_1^2 - a^2 } +
\sqrt {s_2^2
- b^2 } ),}\\
\displaystyle{z_2 = \frac{r} {{s_1 - s_2 }}(\sqrt {s_1^2 - a^2 }
- \sqrt {s_2^2 - b^2 } ),} \\
\displaystyle{w_1 = r\frac{{\sqrt {\Phi (s_2 )} - \sqrt {\Phi (s_1
)} }} {{\Psi
(s_1 ,s_2 ) - \sqrt {\Phi (s_1 )\Phi (s_2 )} }},}\\
\displaystyle{w_2 = r\frac{{\sqrt {\Phi (s_2 )} + \sqrt {\Phi (s_1
)} }} {{\Psi (s_1
,s_2 ) + \sqrt {\Phi (s_1 )\Phi (s_2 )} }},} \\
\displaystyle{w_3 = \frac{1} {{s_1 - s_2 }}[\sqrt {(s_2^2 - b^2
)\Phi (s_1 )} - \sqrt {(s_1^2 - a^2 )\Phi (s_2 )} ].}
\end{array}
$$
From here we immediately obtain explicit expressions for the phase variables
$\alpha_i,\beta_j$, and $\omega_k$ $(i,j,k=1,2,3)$ in terms of the separated variables.

Thus, we have completed the solution of the problem of the
motion of a generalized Kovalevskaya top on the invariant
submanifold $N$.

\vspace{2mm}
R E F E R E N C E S
\vspace{2mm}

1.~\textit{S.\,Kowalevski} Acta Math. {\bf 12}, 177-232 (1889).

2.~\textit{O.\,I.\,Bogoyavlenski\v{\i}}, Dokl. Akad. Nauk SSSR, \textbf{275},
1359–1363 (1984).

3.~\textit{A.\,I.\,Bobenko, A.\,G.\,Reyman, and M.\,A.\,Semenov-Tian-
Shansky}, Commun. Math. Phys., No. 122, 321–354 (1989).

4.~\textit{M.\,P.\,Kharlamov}, Mekh. Tverd. Tela, No. 32, 32–38 (2002).

5.~\textit{D.\,B.\,Zotev}, Regular Chaotic Dyn. \textbf{5} (4), 447–458 (2000).

\vskip5mm

Address: \textit{Russia, 400131, Volgograd, Gagarin Street, 8, VAGS \emph{(}M.P. Kharlamov\emph{)}}

E-mail: \textit{mharlamov@vags.ru}

\flushright{Translated by M.P.Kh.}

\end{document}